\documentclass[11pt]{amsart}
\usepackage[utf8]{inputenc}
\usepackage[T1]{fontenc}
\usepackage{amsmath,amssymb,amsthm}
\usepackage{mathtools}
\usepackage[margin=1.1in]{geometry}
\usepackage{tikz}
\usetikzlibrary{positioning,fit,shapes.geometric,arrows.meta}
\usepackage{enumitem}
\usepackage[colorlinks=true,linkcolor=blue,citecolor=blue,urlcolor=blue]{hyperref}

\newtheorem{theorem}{Theorem}[section]
\newtheorem{lemma}[theorem]{Lemma}
\newtheorem{proposition}[theorem]{Proposition}
\newtheorem{corollary}[theorem]{Corollary}
\newtheorem{conjecture}[theorem]{Conjecture}
\theoremstyle{definition}
\newtheorem{definition}[theorem]{Definition}
\newtheorem{example}[theorem]{Example}
\newtheorem{remark}[theorem]{Remark}

\DeclareMathOperator{\tr}{tr}
\DeclareMathOperator{\diam}{diam}
\newcommand{\eps}{\varepsilon}
\DeclarePairedDelimiter\ceil{\lceil}{\rceil}
\DeclarePairedDelimiter\floor{\lfloor}{\rfloor}
\newcommand{\doi}[1]{\href{https://doi.org/#1}{doi:#1}}

\tikzset{vtx/.style={circle,fill=black,draw,minimum size=5pt,inner sep=0pt},
         wvtx/.style={circle,fill=white,draw,minimum size=5pt,inner sep=0pt}}

\title{Triameter, distance packings, and vertex connectivity}

\author{Artem Hak}
\address{National University of Kyiv-Mohyla Academy, Skovorody str.~2, 04070 Kyiv, Ukraine;\newline
	Kyiv School of Economics, Mykoly Shpaka str.~3, 03113 Kyiv, Ukraine}
\email{artikgak@ukr.net}

\author{Sergiy Kozerenko}
\address{National University of Kyiv-Mohyla Academy, Skovorody str.~2, 04070 Kyiv, Ukraine;\newline
	Kyiv School of Economics, Mykoly Shpaka str.~3, 03113 Kyiv, Ukraine}
\email{kozerenkosergiy@ukr.net}
\thanks{Corresponding author: kozerenkosergiy@ukr.net}


\begin{document}

\begin{abstract}
We study upper bounds for the triameter of a connected graph in terms of distance packings and vertex connectivity. For any positive integer $k$, we bound the triameter of a graph by a linear function of its $k$-packing number and show that this bound is tight for all values of the parameters. As a corollary, we confirm the conjecture of Das and Abdio\u{g}lu on the upper bound for the triameter in terms of the order and the minimum degree of a graph. We also find the maximum triameter among graphs with given order and vertex connectivity, and describe the extremal graphs. Our main tool is a projection of a graph onto a spider or a unicyclic graph which preserves the distances between three prescribed vertices.
\end{abstract}

\maketitle

\section{Introduction}\label{sec-1}

The vertex set of every connected graph $G$ carries the natural
shortest-path metric $d_G$. One of its most basic numerical
characteristics is the diameter
\[
\diam(G)=\max\{d_G(u,v):u,v\in V(G)\}.
\]
A natural three-point analogue of the diameter is obtained by
maximizing the perimeter of a triple of vertices. Namely, the
\emph{triameter} of $G$ is
\[
\tr(G)=
\max\{d_G(u,v)+d_G(u,w)+d_G(v,w):u,v,w\in V(G)\}.
\]
Thus, while the diameter measures how far apart two vertices of a
graph can be, the triameter measures the largest possible total
pairwise distance among three vertices.

The quantity now called the triameter appeared considerably earlier
than its systematic study as a graph invariant. In 1996,
Laurent~\cite{Laurent:96}, in her work on the metric polytope, used the
maximum perimeter of a triple as a normalizing parameter associated
with graph metrics. In particular, she proved the sharp bound
\[
\tr(G)\leq 2n-2
\]
for a connected graph of order $n$, with equality precisely for paths
and trees with three leaves. The same quantity subsequently arose,
apparently independently, in other parts of graph theory. It was used
in lower bounds for radio $k$-chromatic numbers
(see~\cite{KolaPan:15,SahaPan:15}), and Henning and
Yeo~\cite{HenYeo:14} proved that the total domination number satisfies
\[
\gamma_t(G)\geq \frac14
\bigl(d_G(x_1,x_2)+d_G(x_1,x_3)+d_G(x_2,x_3)\bigr)
\]
for every three vertices $x_1,x_2,x_3$, and consequently
$\gamma_t(G)\geq \tr(G)/4$.

A systematic investigation of this parameter in its own right began
only much later. Saha and Basunia~\cite{SahaBas:20} considered the same
three-point metric quantity in connection with radio graceful
labellings, while Das~\cite{Das:21} introduced the term
\emph{triameter} and developed a general theory around it. Among other
results, Das obtained bounds in terms of the order, girth, domination
parameters and vertex connectivity, as well as Nordhaus--Gaddum type
inequalities, and posed several structural and extremal questions.
The order bound $\tr(G)\leq 2n-2$, together with its equality case,
was rediscovered there.

Several of the questions raised by Das were subsequently settled by
Hak, Kozerenko and Oliynyk in~\cite{HakKozOl:22}, where they determined the optimal
lower bound for the triameter of a tree in terms of its order and
number of leaves, and proved two structural properties of triametral
and diametral triples in block graphs. More recently, Das and
Abdio\u{g}lu~\cite{DasAbd:25} characterized graphs with triameter $4$
and $5$, while graphs with small triameters were further investigated
in~\cite{Koz:26}. In particular, the study of small values suggests
that the triameter is sensitive not only to the large-scale metric
structure of a graph, but also to local density and separation
properties.

The present paper is concerned with two such connections: distance
packings and vertex connectivity. Recall that, for a positive integer
$k$, a set $P\subseteq V(G)$ is a \emph{$k$-packing} if
$d_G(u,v)>k$ for every two distinct vertices $u,v\in P$, and
$\rho_k(G)$ denotes the maximum cardinality of a $k$-packing. Distance
packings provide a natural way to measure how many mutually distant
vertices a graph may contain, and therefore it is reasonable to ask
how strongly they constrain the perimeter of a triametral triple.

Our first main result gives a sharp answer to this question.

\begin{theorem}\label{thm-packing}
	Let $G$ be a connected graph with at least three vertices, and let
	$k\in\mathbb{N}$. Then
	\begin{align}\label{eq-packing}
		\tr(G)\leq 2(k+1)\rho_{k}(G)+k-2.
	\end{align}
	Moreover, for all admissible values of the parameters this bound is attained.
\end{theorem}

The case $k=2$ has an immediate consequence for the minimum degree.
Since the closed neighborhoods of the vertices of a $2$-packing are
pairwise disjoint,
\[
\rho_2(G)\leq
\left\lfloor\frac{n}{\delta(G)+1}\right\rfloor .
\]
Hence Theorem~\ref{thm-packing} yields
\[
\tr(G)\leq
6\rho_2(G)
\leq
6\left\lfloor\frac{n}{\delta(G)+1}\right\rfloor .
\]
In particular, this confirms the following conjecture of Das and
Abdio\u{g}lu~\cite{DasAbd:25}.

\begin{conjecture}[Das and Abdio\u{g}lu~\cite{DasAbd:25}]
	\label{conj-DA}
	For any connected graph $G$ with $n$ vertices and minimum degree
	$\delta(G)$,
	\[
	\tr(G)\leq\frac{6n}{\delta(G)+1}.
	\]
\end{conjecture}

The coefficient $6$ in this estimate is asymptotically best possible
for every fixed value of the minimum degree $\delta \ge 2$; see
Proposition~\ref{prop-six}. Theorem~\ref{thm-packing} also produces
bounds in terms of the independence and domination parameters, thus
placing several previously known estimates into a common framework.

Our second direction concerns vertex connectivity. Das~\cite{Das:21}
obtained an upper bound on the triameter in terms of the order and
vertex connectivity. Here we determine the exact extremal value.

\begin{theorem}\label{thm-connectivity}
	Let $n$ and $k$ be positive integers with $n\geq k+2$. Then
	\[
	\max\{\tr(G): |V(G)|=n,\ \kappa(G)=k\}
	=
	2\left\lfloor\frac{n-3}{k}\right\rfloor+4.
	\]
\end{theorem}

We also describe the corresponding extremal graphs. In the
biconnected case this gives the particularly simple bound
$\tr(G)\leq n+1$. Moreover, equality forces a strong structural
restriction: a biconnected graph with $\tr(G)=n+1$ contains a spanning
theta subgraph, and its triametral triple is uniquely determined by
the three branches.

Both main results are based on the same metric construction. Given
three prescribed vertices $a_1,a_2,a_3$ of a graph $G$, we construct
a surjective non-expanding map from $G$ onto a canonical three-legged
model: a spider when the perimeter
$d_G(a_1,a_2)+d_G(a_1,a_3)+d_G(a_2,a_3)$ is even, and a triangle with
three attached paths when it is odd. The map preserves all three
pairwise distances between $a_1,a_2,a_3$. Consequently, distance
packings in the model lift to distance packings in $G$, whereas
suitable vertex cuts of the model pull back to vertex cuts of $G$.
The same two types of graphs were used in~\cite[Proposition~3.1]{Koz:26}
to realize all possible pairs of diameter and triameter. From the
metric viewpoint, the spider model is also closely related to the
tight span, or injective hull, of a three-point metric
space~\cite{BanChe:08,Dress:84,Isbell:64}.

The paper is organized as follows. In Section~\ref{sec-2} we collect
the necessary terminology and preliminary results. Section~\ref{sec-3}
introduces the canonical projection onto a spider or a triangle with
three attached paths and establishes its basic transfer properties.
In Section~\ref{sec-4} we prove Theorem~\ref{thm-packing}; its
consequences for the minimum degree and domination parameters are
derived in Section~\ref{sec-5}. Section~\ref{sec-6} is devoted to
Theorem~\ref{thm-connectivity} and the description of the extremal
graphs. We conclude in Section~\ref{sec-7} with several open
questions.

\section{Definitions and preliminary results}\label{sec-2}

\subsection{General notions}

All graphs in this paper are finite, simple and undirected. For a vertex $u\in V(G)$ in a graph $G$, by $N_{G}(u)=\{x\in V(G):ux\in E(G)\}$ and $N_{G}[u]=N_{G}(u)\cup\{u\}$ we denote the \textit{neighborhood} and the \textit{closed neighborhood} of $u$, respectively. The \textit{degree} of $u$ is the number $\deg_{G}(u)=|N_{G}(u)|$, and $\delta(G)$ denotes the \textit{minimum degree} of $G$. For a set $S\subseteq V(G)$, by $G[S]$ we denote the corresponding induced subgraph, and we put $G-S=G[V(G)\setminus S]$.

A graph is \textit{connected} provided every pair of its vertices is joined by a path. Paths are measured by the number of their edges. In a connected graph $G$, the distance $d_{G}(u,v)$ is the length of a shortest $u-v$ path, and a shortest path is also called a \textit{geodesic}. For a vertex $v\in V(G)$ and an integer $s\geq 0$, we put $B_{G}(v,s)=\{x\in V(G):d_{G}(v,x)\leq s\}$ for the \textit{ball} of radius $s$ around $v$.

A \textit{vertex cut} in a connected graph $G$ is a set $S\subseteq V(G)$ such that $G-S$ is disconnected. A graph $G$ is \textit{$k$-connected} if $|V(G)|>k$ and $G-S$ is connected for every $S\subseteq V(G)$ with $|S|<k$. The \textit{vertex connectivity} $\kappa(G)$ is the largest $k$ such that $G$ is $k$-connected. In particular, every vertex cut in $G$ has at least $\kappa(G)$ vertices. A graph is \textit{biconnected} if it is $2$-connected.

For a triple of vertices $T=\{a_{1},a_{2},a_{3}\}\subseteq V(G)$, a set $S\subseteq V(G)\setminus T$ is called a \textit{$T$-separator} provided the vertices $a_{1},a_{2},a_{3}$ do not all lie in one connected component of $G-S$. Note that a $T$-separator is not required to separate all three vertices $a_{1},a_{2},a_{3}$ from each other. Clearly, every $T$-separator in a connected graph is a vertex cut.

A set $A\subseteq V(G)$ is \textit{independent} if no two of its vertices are adjacent, and $\alpha(G)$ denotes the largest cardinality of an independent set in $G$. As defined in Section~\ref{sec-1}, a set $P\subseteq V(G)$ is a \textit{$k$-packing} if $d_{G}(u,v)>k$ for all distinct $u,v\in P$, and $\rho_{k}(G)$ is the largest cardinality of a $k$-packing. Thus, $\rho_{1}(G)=\alpha(G)$. A set $A\subseteq V(G)$ is \textit{dominating} if every vertex of $G$ lies in $N_{G}[a]$ for some $a\in A$, and $\gamma(G)$ denotes the \textit{domination number} of $G$.

Let $H$ be a graph. A \textit{clique blow-up} of $H$ is a graph obtained from $H$ by replacing every vertex $x\in V(H)$ by a complete graph $K^{(x)}$ with at least one vertex (called the \textit{bag} of $x$), and joining every vertex of $K^{(x)}$ with every vertex of $K^{(y)}$ whenever $xy\in E(H)$. If $H$ is connected, then the distance in the clique blow-up between two vertices from distinct bags $K^{(x)}$, $K^{(y)}$ equals $d_{H}(x,y)$.

\begin{figure}[ht]
\centering
\begin{tikzpicture}[scale=0.9,every node/.style={font=\small}]
\begin{scope}
\node[vtx,label=below:$R$] (c) at (0,0) {};
\draw (c) -- (90:2) node[vtx,label=above:$A_1$] {};
\draw (c) -- (210:2) node[vtx,label=below left:$A_2$] {};
\draw (c) -- (330:2) node[vtx,label=below right:$A_3$] {};
\node at (65:1.1) {$\ell_1$};
\node at (185:1.2) {$\ell_2$};
\node at (355:1.1) {$\ell_3$};
\node at (0,-2.3) {$S(\ell_1,\ell_2,\ell_3)$};
\end{scope}
\begin{scope}[xshift=6.5cm]
\coordinate (R1) at (90:0.5); \coordinate (R2) at (210:0.5); \coordinate (R3) at (330:0.5);
\draw (R1)--(R2)--(R3)--cycle;
\foreach \p in {R1,R2,R3} \node[vtx] at (\p) {};
\node[right] at (R1) {$R_1$};
\node[above left] at (R2) {$R_2$};
\node[above right] at (R3) {$R_3$};
\draw (R1) -- (90:2.3) node[vtx,label=above:$A_1$] {};
\draw (R2) -- (210:2.3) node[vtx,label=below left:$A_2$] {};
\draw (R3) -- (330:2.3) node[vtx,label=below right:$A_3$] {};
\node at (75:1.5) {$\ell_1$};
\node at (195:1.6) {$\ell_2$};
\node at (345:1.6) {$\ell_3$};
\node at (0,-2.3) {$U(\ell_1,\ell_2,\ell_3)$};
\end{scope}
\end{tikzpicture}
\caption{The spider $S(\ell_1,\ell_2,\ell_3)$ and the graph $U(\ell_1,\ell_2,\ell_3)$. The labeled strokes are paths of the indicated lengths.}
\label{fig-spiders}
\end{figure}
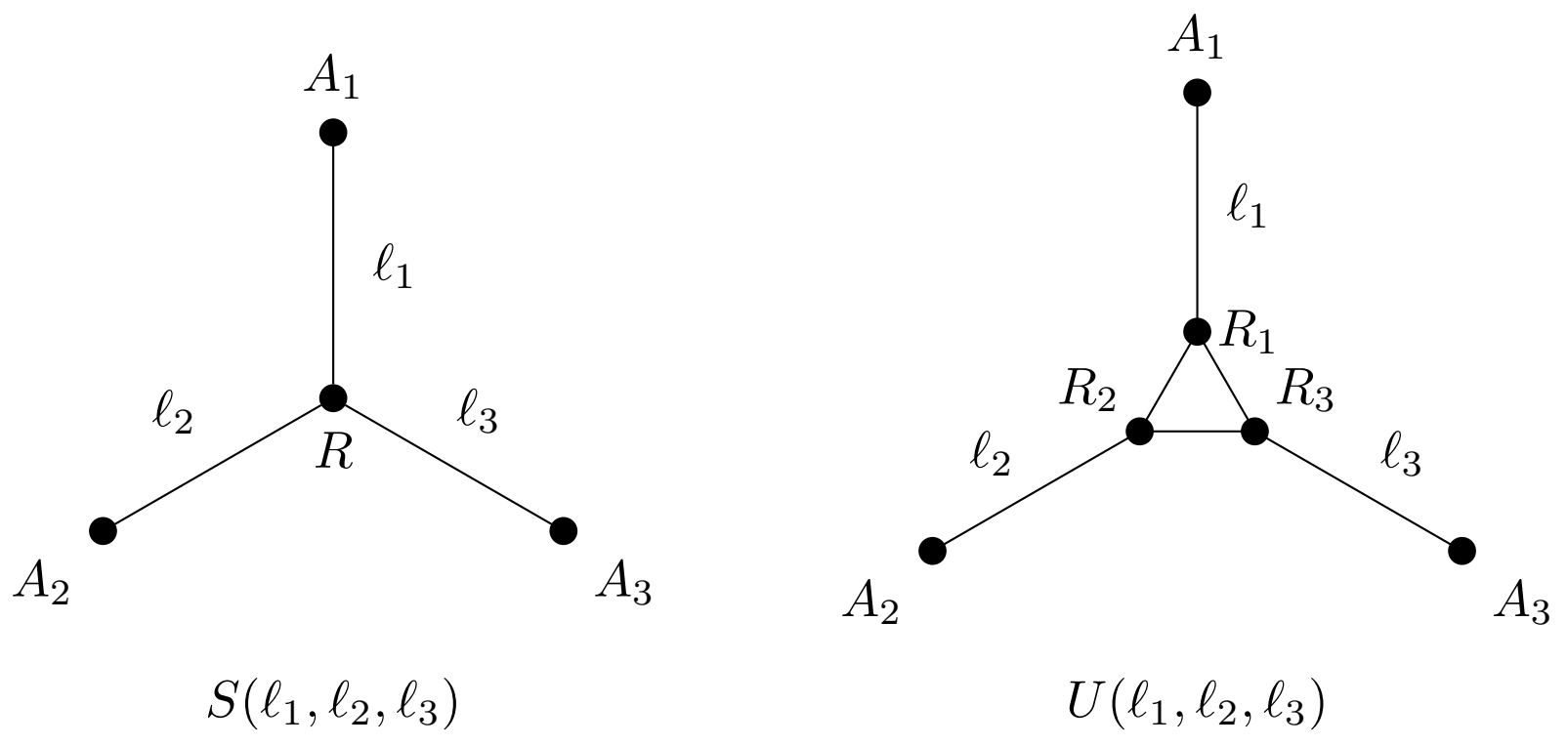

Let $\ell_{1},\ell_{2},\ell_{3}\geq 0$ be integers. The \textit{spider} $S(\ell_{1},\ell_{2},\ell_{3})$ is obtained from three paths of lengths $\ell_{1},\ell_{2},\ell_{3}$ by identifying one end of each path into a single vertex $R$, called the \textit{center}. The three paths are called the \textit{legs} of the spider, and the other ends $A_{1},A_{2},A_{3}$ of the legs are called its \textit{feet} (if $\ell_{i}=0$, then $A_{i}=R$). Similarly, the graph $U(\ell_{1},\ell_{2},\ell_{3})$ is obtained from the triangle $R_{1}R_{2}R_{3}$ by attaching a path of length $\ell_{i}$ at the vertex $R_{i}$ for $i=1,2,3$. The path attached at $R_{i}$ together with $R_{i}$ is called the $i$--th \textit{leg} of $U(\ell_{1},\ell_{2},\ell_{3})$, and the end $A_{i}$ of this path is its $i$--th \textit{foot} (again, $A_{i}=R_{i}$ if $\ell_{i}=0$). In both cases, the \textit{depth} of a vertex on the $i$--th leg is its distance to $R$ (for spiders) or to $R_{i}$ (for the graphs $U(\ell_{1},\ell_{2},\ell_{3})$). Thus, the foot $A_{i}$ has depth $\ell_{i}$, and the center of a spider has depth $0$ on each of its legs. Both types of graphs are depicted in Figure~\ref{fig-spiders}. Clearly,
\begin{align}\label{eq-order-spider}
|V(S(\ell_{1},\ell_{2},\ell_{3}))|=\ell_{1}+\ell_{2}+\ell_{3}+1,\quad |V(U(\ell_{1},\ell_{2},\ell_{3}))|=\ell_{1}+\ell_{2}+\ell_{3}+3.
\end{align}

A \textit{theta graph} consists of two vertices (called \textit{hubs}) joined by three internally disjoint paths of lengths at least two. A \textit{subdivision} of a graph $H$ is obtained by replacing the edges of $H$ by internally disjoint paths.

Finally, let $G$ and $H$ be connected graphs. A map $\varphi\colon V(G)\to V(H)$ is called \textit{non-expanding} if $d_{H}(\varphi(u),\varphi(v))\leq d_{G}(u,v)$ for all $u,v\in V(G)$. Such maps are also known as \textit{metric maps} or \textit{$1$-Lipschitz maps}; we will use the term non-expanding throughout the paper. One can show that a map $\varphi$ is non-expanding if and only if for every edge $uv\in E(G)$, either $\varphi(u)=\varphi(v)$ or $\varphi(u)\varphi(v)\in E(H)$.

\subsection{Triameter}

For a triple of vertices $u,v,w\in V(G)$ in a connected graph $G$, we put
\[d_{G}(u,v,w)=d_{G}(u,v)+d_{G}(u,w)+d_{G}(v,w),\]
and call this number the \textit{perimeter} of the triple. The \textit{triameter} of $G$ is the number
\[\tr(G)=\max\{d_{G}(u,v,w):u,v,w\in V(G)\}.\]
A triple of distinct vertices $u,v,w$ with $d_{G}(u,v,w)=\tr(G)$ is called \textit{triametral}. Note that allowing repeated vertices does not change the maximum: if $|V(G)|\geq 3$, then $d_{G}(u,u,v)=2d_{G}(u,v)\leq d_{G}(u,v,w)$ for any third vertex $w$. The following bounds are well known.

\begin{proposition}{\rm\bf(folklore)}\label{prop-trivial}
For any connected graph $G$, we have $2\diam(G)\leq\tr(G)\leq 3\diam(G)$.
\end{proposition}

\begin{proposition}{\rm\cite{Laurent:96,Das:21}}\label{prop-Das-order}
For any connected graph $G$ with $n\geq 3$ vertices, we have $\tr(G)\leq 2n-2$. Moreover, the equality holds if and only if $G$ is a path or a tree with exactly three leaves.
\end{proposition}

\begin{proposition}{\rm\cite[Corollary~3.11]{Das:21}}\label{prop-Das-dom}
For any connected graph $G$, we have $\tr(G)\leq 6\gamma(G)$.
\end{proposition}

\begin{proposition}{\rm\cite[Proposition~3.14]{Das:21}}\label{prop-Das-kappa}
For any connected graph $G$ with $n$ vertices, we have $\tr(G)\leq\frac{3(n-2)}{\kappa(G)}+3$.
\end{proposition}

We will need the triameters of spiders and of the graphs $U(\ell_{1},\ell_{2},\ell_{3})$.

\begin{lemma}\label{lem-tr-spider}
Let $\ell_{1},\ell_{2},\ell_{3}\geq 0$ be integers, and let $A_{1},A_{2},A_{3}$ be the feet. Then
\begin{enumerate}
\item $d(A_{i},A_{j})=\ell_{i}+\ell_{j}$ and $\tr(S(\ell_{1},\ell_{2},\ell_{3}))=d(A_{1},A_{2},A_{3})=2(\ell_{1}+\ell_{2}+\ell_{3})$, where $d$ is the distance in $S(\ell_{1},\ell_{2},\ell_{3})$;
\item $d(A_{i},A_{j})=\ell_{i}+\ell_{j}+1$ and $\tr(U(\ell_{1},\ell_{2},\ell_{3}))=d(A_{1},A_{2},A_{3})=2(\ell_{1}+\ell_{2}+\ell_{3})+3$, where $d$ is the distance in $U(\ell_{1},\ell_{2},\ell_{3})$.
\end{enumerate}
Moreover, in both graphs the geodesic between any two feet is unique.
\end{lemma}
\begin{proof}
The formulas for $d(A_{i},A_{j})$ and the uniqueness of geodesics are clear: in $U(\ell_{1},\ell_{2},\ell_{3})$, the route from $A_{i}$ to $A_{j}$ through the third triangle vertex $R_{k}$ is longer by one than the route through the edge $R_{i}R_{j}$.

1. Let $S=S(\ell_{1},\ell_{2},\ell_{3})$ and fix $x,y,z\in V(S)$. Since $S$ is a tree, the distance between two vertices equals the number of edges whose deletion separates them. Each edge $e\in E(S)$ separates at most two of the three pairs $\{x,y\},\{x,z\},\{y,z\}$. Hence, $d_{S}(x,y,z)\leq 2|E(S)|=2(\ell_{1}+\ell_{2}+\ell_{3})$, and the feet attain this value.

2. Let $U=U(\ell_{1},\ell_{2},\ell_{3})$ and fix $x,y,z\in V(U)$. For two vertices $u,v\in V(U)$, the distance $d_{U}(u,v)$ equals the number of leg edges separating $u$ and $v$, plus one if $u$ and $v$ lie on different legs. As before, each leg edge contributes at most $2$ to the perimeter $d_U(x,y,z)$. In addition, for each of the three
pairs, a shortest path may contain at most one edge of the central triangle. Therefore, $d_{U}(x,y,z)\leq 2(\ell_{1}+\ell_{2}+\ell_{3})+3$. Equality is attained for the triple $(A_1,A_2,A_3)$.
\end{proof}

\section{Projection onto a spider}\label{sec-3}

Throughout this section, $G$ is a connected graph and
$a_1,a_2,a_3\in V(G)$ are three distinct vertices. Put
$T=\{a_1,a_2,a_3\}$ and
\[
p=d_G(a_1,a_2)+d_G(a_1,a_3)+d_G(a_2,a_3).
\]
Our goal is to associate with $T$ a canonical graph $H_T$ that preserves
the three pairwise distances between $a_1,a_2,a_3$ and is as small as
possible. Depending on the parity of $p$, this graph will be either a
spider or a triangle with three attached paths.

Let $\eps\in\{0,1\}$ be the remainder of $p$ modulo $2$. For
$\{i,j,k\}=\{1,2,3\}$ define
\[
\ell_i=
\frac{
	d_G(a_i,a_j)+d_G(a_i,a_k)-d_G(a_j,a_k)-\eps
}{2}.
\]
The numbers $\ell_1,\ell_2,\ell_3$ will be the lengths of the three
legs of $H_T$.

By the triangle inequality,
$d_G(a_i,a_j)+d_G(a_i,a_k)-d_G(a_j,a_k)$ is nonnegative, and it has
the same parity as $p$. Moreover, if it is zero, then $p$ is even and
$\eps=0$. Hence $\ell_1,\ell_2,\ell_3$ are nonnegative integers.
A direct computation gives
\begin{equation}\label{eq-ell}
\begin{aligned}
	d_G(a_i,a_j)&=\ell_i+\ell_j+\eps
	\quad\text{for } i\neq j,\\
	\ell_1+\ell_2+\ell_3&=\frac{p-3\eps}{2}.
\end{aligned}
\end{equation}
If $\eps=0$, then at most one of $\ell_1,\ell_2,\ell_3$ is zero;
otherwise $d_G(a_i,a_j)=0$ for some $i\neq j$.

\begin{definition}\label{def-model}
	The \textit{canonical model} of the triple $T$ is the graph
	\[
	H_T=
	\begin{cases}
		S(\ell_1,\ell_2,\ell_3), & \text{if } p \text{ is even},\\
		U(\ell_1,\ell_2,\ell_3), & \text{if } p \text{ is odd},
	\end{cases}
	\]
	with distinguished feet $A_1,A_2,A_3$.
\end{definition}

By Lemma~\ref{lem-tr-spider} and~\eqref{eq-ell},
\[
d_{H_T}(A_i,A_j)=d_G(a_i,a_j)
\qquad\text{for all }i\neq j.
\]
Thus $H_T$ realizes exactly the three pairwise distances determined by
$T$. Moreover,
\[
|V(H_T)|=\ceil*{\frac{p}{2}}+1.
\]
If $p$ is even and one of the $\ell_i$ is zero, then $H_T$ is a path
and the corresponding foot $A_i$ is an interior vertex.

The following lemma shows that this canonical model is not only
distance-preserving on the distinguished triple, but also arises as a
non-expanding quotient of the whole graph $G$.

\begin{lemma}\label{lem-projection}
	There exists a surjective non-expanding map
	$\varphi\colon V(G)\to V(H_T)$ with $\varphi(a_i)=A_i$ for
	$i=1,2,3$. Moreover, every edge of $H_T$ is of the form
	$\varphi(u)\varphi(v)$ for some edge $uv\in E(G)$, $|V(H_T)|=\ceil*{\frac{p}{2}}+1$, and
	\[
	\tr(H_T)=p.
	\]
\end{lemma}
\begin{proof}
	We write $H=H_{T}$ and construct the map $\varphi$ separately for even and odd $p$.
	
	\textbf{Case 1:} $p$ is even. We first observe that, for every $v\in V(G)$, the inequality $d_G(a_i,v)<\ell_i$ can hold for at most one index $i$. Indeed, if
	$d_G(a_i,v)<\ell_i$ and $d_G(a_j,v)<\ell_j$ for some $i\neq j$, then
	\[
	d_G(a_i,a_j)\leq d_G(a_i,v)+d_G(v,a_j)
	<\ell_i+\ell_j=d_G(a_i,a_j),
	\]
	a contradiction. Thus, if $d_G(a_i,v)<\ell_i$ for some $i$, let
	$\varphi(v)$ be the vertex of depth $\ell_i-d_G(a_i,v)$ on the
	$i$--th leg of $H$. If no such index exists, let $\varphi(v)=R$, the
	center of $H$.
	
	\textbf{Case 2:} $p$ is odd. We claim that the balls $B_{G}(a_{1},\ell_{1})$, $B_{G}(a_{2},\ell_{2})$, $B_{G}(a_{3},\ell_{3})$ are pairwise disjoint. Indeed, a vertex $v\in B_{G}(a_{i},\ell_{i})\cap B_{G}(a_{j},\ell_{j})$ for $i\neq j$ would give
	\[d_{G}(a_{i},a_{j})\leq d_{G}(a_{i},v)+d_{G}(v,a_{j})\leq\ell_{i}+\ell_{j}<d_{G}(a_{i},a_{j}).\]
	Now, for $v\in B_{G}(a_{i},\ell_{i})$, let $\varphi(v)$ be the vertex of depth $\ell_{i}-d_{G}(a_{i},v)$ on the $i$--th leg of $H$ (in particular, the depth $0$ corresponds to $R_{i}$). Every vertex outside the three balls is mapped to $R_{1}$.
	
	In both cases $\varphi(a_{i})$ is the vertex of depth $\ell_{i}$ on the $i$--th leg, that is, $\varphi(a_{i})=A_{i}$. We split the rest of the proof into a series of claims.
	
	\textbf{Claim 1.} If $\varphi(w)$ lies on the $i$--th leg of $H$, then its depth equals $\max\{\ell_{i}-d_{G}(a_{i},w),0\}$.
	
	For the vertices $w$ with $\ell_{i}-d_{G}(a_{i},w)\geq 0$ this is the definition of $\varphi$. In Case~1, a vertex mapped to the center satisfies $\ell_{i}-d_{G}(a_{i},w)\leq 0$ for all $i$, and the center has depth $0$ on every leg. In Case~2, a vertex outside the three balls is mapped to $R_{1}$, which has depth $0$ on the first leg, and it satisfies $\ell_{1}-d_{G}(a_{1},w)<0$. This proves the claim.
	
	\textbf{Claim 2.} The map $\varphi$ is non-expanding.
	
	Fix $u,v\in V(G)$. At first, suppose that $\varphi(u)$ and $\varphi(v)$ lie on a common leg $i$. Then $d_{H}(\varphi(u),\varphi(v))$ equals the difference of their depths. Since the function $x\mapsto\max\{x,0\}$ does not increase differences, Claim~1 and the triangle inequality give
	\[d_{H}(\varphi(u),\varphi(v))\leq|d_{G}(a_{i},u)-d_{G}(a_{i},v)|\leq d_{G}(u,v).\]
	Now suppose that $\varphi(u)$ and $\varphi(v)$ lie on distinct legs of
	$H$, say the $i$--th and $j$--th legs, respectively. Let $x$ and $y$
	denote their respective depths.
	
	If $p$ is even, then $H$ is a spider. Since a vertex of depth $0$ is
	the center of $H$, which belongs to every leg, the assumption that
	$\varphi(u)$ and $\varphi(v)$ do not lie on a common leg implies that
	$x,y>0$. Hence $x=\ell_i-d_G(a_i,u)$, $y=\ell_j-d_G(a_j,v)$, and
	$d_H(\varphi(u),\varphi(v))=x+y$. By the triangle inequality,
	\begin{align*}
		\ell_{i}+\ell_{j}&=d_{G}(a_{i},a_{j})\leq d_{G}(a_{i},u)+d_{G}(u,v)+d_{G}(v,a_{j})\\
		&=(\ell_{i}-x)+d_{G}(u,v)+(\ell_{j}-y),
	\end{align*}
	which implies $d_{G}(u,v)\geq x+y$.
	
	If $p$ is odd, then $d_H(\varphi(u),\varphi(v))=x+y+1$. If $u\in B_G(a_i,\ell_i)$ and $v\in B_G(a_j,\ell_j)$, then the same
	computation, using
	\[
	d_G(a_i,a_j)=\ell_i+\ell_j+1,
	\]
	gives
	\[
	d_G(u,v)\geq x+y+1.
	\]
	
	Otherwise, suppose that $u$ lies outside the three balls. By the
	construction of $\varphi$, we then have $\varphi(u)=R_1$, so that
	$i=1$ and $x=0$. Since $\varphi(u)\neq\varphi(v)$, it follows that
	$j\neq 1$, and hence $v\in B_G(a_j,\ell_j)$. Therefore,
	\[
	d_G(u,v)
	\geq d_G(a_j,u)-d_G(a_j,v)
	\geq (\ell_j+1)-(\ell_j-y)
	= y+1
	= d_H(\varphi(u),\varphi(v)).
	\]
	
	Thus, in all cases, $d_H(\varphi(u),\varphi(v))\leq d_G(u,v)$, and the claim follows.
	
	\textbf{Claim 3.} The map $\varphi$ is surjective, and every edge of $H$ is of the form $\varphi(u)\varphi(v)$ for some $uv\in E(G)$.
	
	Fix $i\neq j$ and a geodesic $a_{i}=x_{0},x_{1},\dots,x_{m}=a_{j}$ in $G$, where $m=d_{G}(a_{i},a_{j})$. By Claim~2, consecutive vertices $\varphi(x_{r}),\varphi(x_{r+1})$ are equal or adjacent in $H$. Thus, $\varphi(x_{0}),\dots,\varphi(x_{m})$ is a walk in $H$ from $A_{i}$ to $A_{j}$ of length $m=d_{H}(A_{i},A_{j})$. Therefore, this walk is a geodesic in $H$, and all its steps are edges. By Lemma~\ref{lem-tr-spider}, it is the unique geodesic between $A_{i}$ and $A_{j}$. It remains to note that the three geodesics between the feet of $H$ cover all vertices and all edges of $H$. This proves the claim.
	
	Finally, by~\eqref{eq-order-spider} and~\eqref{eq-ell}, $|V(H)|=\frac{p}{2}+1$ if $p$ is even, and $|V(H)|=\frac{p-3}{2}+3=\frac{p+3}{2}$ if $p$ is odd. In both cases $|V(H)|=\ceil*{\frac{p}{2}}+1$. The equality $\tr(H)=p$ follows from Lemma~\ref{lem-tr-spider} and~\eqref{eq-ell}.
\end{proof}

The map $\varphi$ from Lemma~\ref{lem-projection} will be called the \textit{projection} of $G$ onto $H_{T}$. Note that the preimages $\varphi^{-1}(x)$ of vertices $x\in V(H_{T})$ need not induce connected subgraphs of $G$.

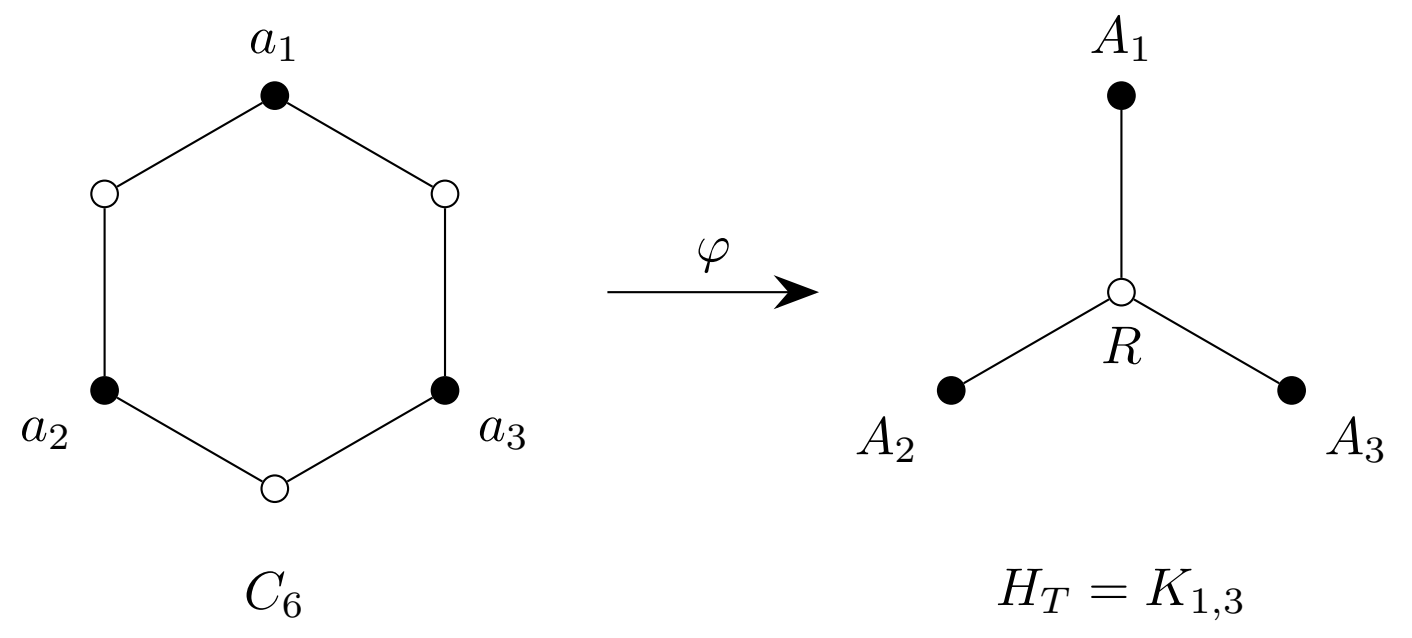
\begin{figure}[ht]
\centering
\begin{tikzpicture}[scale=1,every node/.style={font=\small}]
\begin{scope}
\foreach \a/\n in {90/1,210/2,330/3} \node[vtx,label={\a:$a_{\n}$}] (a\n) at (\a:1.3) {};
\foreach \a/\n in {30/1,150/2,270/3} \node[wvtx] (b\n) at (\a:1.3) {};
\draw (a1)--(b1)--(a3)--(b3)--(a2)--(b2)--(a1);
\node at (0,-2) {$C_6$};
\end{scope}
\draw[-{Stealth[length=3mm]}] (2.2,0) -- node[above] {$\varphi$} (3.6,0);
\begin{scope}[xshift=5.6cm]
\node[wvtx,label=below:$R$] (c) at (0,0) {};
\foreach \a/\n in {90/1,210/2,330/3} {\node[vtx,label={\a:$A_{\n}$}] (A\n) at (\a:1.3) {}; \draw (c)--(A\n);}
\node at (0,-2) {$H_T=K_{1,3}$};
\end{scope}
\end{tikzpicture}
\caption{The projection from Example~\ref{ex-C6}: the three white vertices of $C_6$ are mapped to the center $R$.}
\label{fig-C6}
\end{figure}

\begin{example}\label{ex-C6}
Let $a_{1},a_{2},a_{3}$ be three pairwise nonadjacent vertices of the cycle $C_{6}$. Then $p=6$ and $\ell_{1}=\ell_{2}=\ell_{3}=1$, so $H_{T}=S(1,1,1)=K_{1,3}$. The preimage of the center consists of the other three vertices of $C_{6}$, which are also pairwise nonadjacent (see Figure~\ref{fig-C6}).
\end{example}

The following lemma explains why the projection is useful.

\begin{lemma}\label{lem-transfer}
	Let $\varphi\colon V(G)\to V(H_T)$ be the projection from
	Lemma~\ref{lem-projection}. Then the following hold.
	\begin{enumerate}
		\item If $P$ is a $k$-packing in $H_T$, then choosing one preimage of
		each vertex of $P$ yields a $k$-packing in $G$ of the same cardinality.
		In particular,
		\[
		\rho_k(H_T)\leq \rho_k(G).
		\]
		
		\item If $X$ is a $T$-separator in $H_T$, then $\varphi^{-1}(X)$ is a $T$-separator in $G$.
	\end{enumerate}
\end{lemma}
\begin{proof}
1. If $u,v\in V(G)$ are preimages of two distinct vertices of $P$, then
\[d_{G}(u,v)\geq d_{H_{T}}(\varphi(u),\varphi(v))>k.\]

2. Since $\varphi(a_{i})=A_{i}\notin X$, the set $\varphi^{-1}(X)$ does not contain $a_{1},a_{2},a_{3}$. It is nonempty since $\varphi$ is surjective. Now suppose that $a_{i}$ and $a_{j}$ are joined by a path $a_{i}=x_{0},x_{1},\dots,x_{m}=a_{j}$ in $G-\varphi^{-1}(X)$. Since $\varphi$ is non-expanding, the sequence $\varphi(x_{0}),\dots,\varphi(x_{m})$ is a walk from $A_{i}$ to $A_{j}$ in $H_{T}-X$. Therefore, if $A_{i}$ and $A_{j}$ lie in different connected components of $H_{T}-X$, then so do $a_{i}$ and $a_{j}$ in $G-\varphi^{-1}(X)$.
\end{proof}

As a first application, we give an alternative proof of Proposition~\ref{prop-Das-order}.

\begin{proof}[Alternative proof of Proposition~\ref{prop-Das-order}]
Let $a_{1},a_{2},a_{3}$ be a triametral triple in $G$, and let $\varphi\colon V(G)\to V(H_{T})$ be the projection. Since $\varphi$ is surjective, Lemma~\ref{lem-projection} gives $\ceil*{\frac{\tr(G)}{2}}+1\leq n$, which implies $\tr(G)\leq 2n-2$.

Now suppose that $\tr(G)=2n-2$. Then $p=\tr(G)$ is even and $|V(H_{T})|=n$, so $\varphi$ is a bijection. Since $\varphi$ is non-expanding and injective, it maps every edge of $G$ to an edge of $H_{T}$. Conversely, by Lemma~\ref{lem-projection}, every edge of $H_{T}$ is the image of an edge of $G$. Hence, $\varphi$ is an isomorphism, and $G\simeq H_{T}$ is a path or a tree with three leaves. Conversely, for such trees the feet of the corresponding spider (or the two ends of the path and any third vertex) form a triple with perimeter $2(n-1)$.
\end{proof}

\section{Distance packings}\label{sec-4}

In this section, we prove Theorem~\ref{thm-packing}. We keep the notation of Section~\ref{sec-3}. First, we find a large $k$-packing in the canonical model $H_T$ associated with a given triple $T$. We then use this packing, together with the distance-preserving projection constructed in Section~\ref{sec-3}, to prove Theorem~\ref{thm-packing}.

\begin{lemma}\label{lem-packing-spider}
	Let $k\in\mathbb{N}$. Write $\ell_{i}=(k+1)q_{i}+s_{i}$ with $q_{i}\geq 0$ and $0\leq s_{i}\leq k$ for $i=1,2,3$, and put $q=q_{1}+q_{2}+q_{3}$. Then
	\begin{align}\label{eq-packing-spider}
		p\leq 2(k+1)\rho_{k}(H_{T})+k-2.
	\end{align}
	Moreover, if the equality holds in~\eqref{eq-packing-spider}, then $\rho_{k}(H_{T})=q+1$ and $d_{G}(a_{i},a_{j})=(k+1)(q_{i}+q_{j})+k$ for all $i\neq j$.
\end{lemma}
\begin{proof}
	Put $H=H_T$. For $i=1,2,3$, let $u_i$ be the vertex of depth $s_i$
	on the $i$--th leg of $H$; if $p$ is even and $s_i=0$, then $u_i$ is
	the center. On the $i$--th leg choose the $q_i$ vertices of depths
	\[
	s_i+(k+1),\ s_i+2(k+1),\ \dots,\ s_i+q_i(k+1),
	\]
	and let $Q$ be the set of all chosen vertices. Thus $|Q|=q$.
	
	\textbf{Claim 1.} The set $Q$ is a $k$-packing in $H$, and
	$d_H(x,u_j)>k$ for every $x\in Q$ and $j\in\{1,2,3\}$.
	
	Indeed, two chosen vertices on the same leg are at distance at least
	$k+1$. Two chosen vertices on different legs have depths at least
	$k+1$, and hence their distance is greater than $k$. Moreover, if
	$x\in Q$ lies on the $i$--th leg, then
	$d_H(x,u_i)\geq k+1$, while for $j\neq i$ the depth of $x$ alone is
	already greater than $k$, so $d_H(x,u_j)>k$. This proves the claim.
	
	Thus, every $u_i$ can be added individually to $Q$. However, the three
	vertices $u_1,u_2,u_3$ need not form a $k$-packing, since some of their
	mutual distances may be at most $k$. We therefore have to determine
	how many of them can be added simultaneously.
	
	Put $p_0=d_H(u_1,u_2,u_3)$. Since
	$d_H(u_i,u_j)=s_i+s_j+\eps$ for $i\neq j$, we have
	\[
	p_0=2(s_1+s_2+s_3)+3\eps.
	\]
	On the other hand, substituting
	$\ell_i=(k+1)q_i+s_i$ into the second identity in~\eqref{eq-ell} gives
	\begin{align*}
		p
		&=2(\ell_1+\ell_2+\ell_3)+3\eps\\
		&=2(k+1)(q_1+q_2+q_3)+2(s_1+s_2+s_3)+3\eps\\
		&=2(k+1)q+p_0.
	\end{align*}
	
	Let $b\in\{1,2,3\}$ be the maximum cardinality of a $k$-packing
	contained in $\{u_1,u_2,u_3\}$. Our next goal is to bound $p_0$ in
	terms of $b$ in such a way that the term $2(k+1)b$ can be combined
	with $2(k+1)q$ in the above expression for $p$.
	
	\textbf{Claim 2.} We have
	\[
	p_0\leq 2(k+1)b+k-2,
	\]
	and the inequality is strict whenever $b\geq2$.
	
	We consider the three possible values of $b$.
	
	\textbf{Case 1:} $b=1$. No two of $u_1,u_2,u_3$ can belong to a
	$k$-packing simultaneously. Hence
	$d_H(u_i,u_j)\leq k$ for every $i\neq j$, and therefore
	\[
	p_0=d_H(u_1,u_2)+d_H(u_1,u_3)+d_H(u_2,u_3)
	\leq3k.
	\]
	Substituting $b=1$ into the claimed bound gives
	$2(k+1)+k-2=3k$, as required.
	
	\textbf{Case 2:} $b=2$. Since all three vertices cannot belong to a
	$k$-packing, some pair, say $u_i,u_j$, satisfies
	$d_H(u_i,u_j)\leq k$. Thus
	$s_i+s_j+\eps\leq k$. If $l$ is the remaining index, then
	$s_l\leq k$, and hence
	\begin{align*}
		p_0
		&=2(s_i+s_j+\eps)+2s_l+\eps\\
		&\leq4k+1.
	\end{align*}
	For $b=2$, the right-hand side in Claim~2 is
	$2(k+1)\cdot2+k-2=5k+2$. Hence
	$p_0<5k+2$, so the inequality is strict.
	
	\textbf{Case 3:} $b=3$. Since $s_i\leq k$ for each $i$, we obtain
	$p_0=2(s_1+s_2+s_3)+3\eps\leq6k+3$. On the other hand, substituting
	$b=3$ into the right-hand side of Claim~2 gives
	$2(k+1)\cdot3+k-2=7k+4$. Thus
	$p_0<7k+4$, and the inequality is again strict.
	
	This proves Claim~2.
	
	By Claim~1, a largest $k$-packing contained in
	$\{u_1,u_2,u_3\}$ can be added to $Q$. Hence $H$ contains a
	$k$-packing $P$ with $|P|=q+b$, and therefore
	$\rho_k(H)\geq q+b$. Combining the identity
	$p=2(k+1)q+p_0$ with Claim~2, we obtain
	\begin{align*}
		p
		&=2(k+1)q+p_0\\
		&\leq2(k+1)(q+b)+k-2\\
		&\leq2(k+1)\rho_k(H)+k-2.
	\end{align*}
	This proves~\eqref{eq-packing-spider}.
	
	Finally, suppose that equality holds in~\eqref{eq-packing-spider}.
	Then equality must hold in both inequalities above. In particular,
	equality holds in Claim~2. Since the inequality in Claim~2 is strict
	for $b=2$ and $b=3$, we must have $b=1$. Consequently,
	$p_0=3k$, and equality in the second inequality gives
	$\rho_k(H)=q+b=q+1$.
	
	Since $b=1$, every pair $u_i,u_j$ has distance at most $k$. Their
	three pairwise distances sum to $p_0=3k$, so in fact
	$d_H(u_i,u_j)=k$ for every $i\neq j$. Hence
	$s_i+s_j+\eps=k$ for all $i\neq j$. Finally, using the first identity
	in~\eqref{eq-ell}, we obtain
	\begin{align*}
		d_G(a_i,a_j)
		&=\ell_i+\ell_j+\eps\\
		&=(k+1)(q_i+q_j)+s_i+s_j+\eps\\
		&=(k+1)(q_i+q_j)+k
	\end{align*}
	for all $i\neq j$.
\end{proof}

Now we are ready to prove Theorem~\ref{thm-packing}.

\begin{proof}[Proof of Theorem~\ref{thm-packing}]
	Let $a_1,a_2,a_3$ be a triametral triple in $G$. Thus, for
	$T=\{a_1,a_2,a_3\}$, we have
	$p=d_G(a_1,a_2,a_3)=\tr(G)$. By Lemma~\ref{lem-packing-spider} and
	Lemma~\ref{lem-transfer}(1),
	\[
	\tr(G)=p
	\leq 2(k+1)\rho_k(H_T)+k-2
	\leq 2(k+1)\rho_k(G)+k-2.
	\]
	This proves the required bound.
	
	Suppose now that equality holds in~\eqref{eq-packing}. Then both
	inequalities above are equalities. Hence
	$\rho_k(H_T)=\rho_k(G)$, and equality holds in
	\eqref{eq-packing-spider}. By Lemma~\ref{lem-packing-spider}, we have
	$\rho_k(H_T)=q+1$, where $q=q_1+q_2+q_3$, and
	\[
	d_G(a_i,a_j)=(k+1)(q_i+q_j)+k
	\]
	for all $i\neq j$. Therefore,
	\[
	q_1+q_2+q_3=\rho_k(G)-1.
	\]
	
	It remains to show that the bound is sharp for all admissible values of
	$k$ and $\rho_k(G)$. Fix $k,h\in\mathbb{N}$ and choose nonnegative
	integers $q_1,q_2,q_3$ such that
	$q_1+q_2+q_3=h-1$. Write
	$k=2s+\eps$, where $s=\floor*{k/2}$ and $\eps\in\{0,1\}$.
	
	We construct $G_{k,h}$ by taking
	\[
	G_{k,h}=
	\begin{cases}
		S\bigl((k+1)q_1+s,\,(k+1)q_2+s,\,(k+1)q_3+s\bigr),
		& \text{if } \eps=0,\\
		U\bigl((k+1)q_1+s,\,(k+1)q_2+s,\,(k+1)q_3+s\bigr),
		& \text{if } \eps=1.
	\end{cases}
	\]
	Let $A_1,A_2,A_3$ be its feet. By Lemma~\ref{lem-tr-spider},
	\begin{align*}
		\tr(G_{k,h})
		&=2\sum_{i=1}^{3}\bigl((k+1)q_i+s\bigr)+3\eps\\
		&=2(k+1)(h-1)+6s+3\eps\\
		&=2(k+1)h+k-2.
	\end{align*}
	
	We next show that $\rho_k(G_{k,h})=h$. On each leg, separate the first
	$s$ vertices from the remaining $(k+1)q_i$ vertices, and divide the
	latter into $q_i$ consecutive blocks of $k+1$ vertices each. Together
	with the center $R$ when $\eps=0$, or the triangle
	$R_1R_2R_3$ when $\eps=1$, the initial segments of the three legs form
	a central part.
	
	The central part has diameter $2s+\eps=k$, while every block has
	diameter $k$. Consequently, a $k$-packing contains at most one vertex
	from each block and at most one vertex from the central part. Since
	there are
	\[
	1+q_1+q_2+q_3=h
	\]
	such parts, we obtain $\rho_k(G_{k,h})\leq h$.
	
	For the reverse inequality, choose one vertex from the central part,
	namely $R$ if $\eps=0$ and $R_1$ if $\eps=1$, and on each leg choose
	the deepest vertex of every block. This gives $h$ vertices in total.
	Two chosen vertices from consecutive blocks on the same leg are at
	distance $k+1$. The chosen vertex in the first block of the $i$--th
	leg has depth $s+k+1$, and therefore its distance from the chosen
	central vertex is greater than $k$. Finally, chosen vertices lying on
	different legs are also at distance greater than $k$. Hence these
	vertices form a $k$-packing, so $\rho_k(G_{k,h})\geq h$.
	
	Therefore,
	\[
	\rho_k(G_{k,h})=h
	\]
	and
	\[
	\tr(G_{k,h})
	=2(k+1)\rho_k(G_{k,h})+k-2.
	\]
	Thus the bound in~\eqref{eq-packing} is sharp.
\end{proof}

\begin{figure}[ht]
	\centering
	\begin{tikzpicture}[scale=0.9,every node/.style={font=\small}]
		
		\node[wvtx,label={[label distance=8pt]below:$R$}] (c) at (0,0) {};
		
		\foreach \i in {1,2,3,5,6} {
			\node[vtx] (l\i) at (-\i,0) {};
		}
        
        \node[wvtx] (l4) at (-4,0) {};
        
		\node[wvtx] (l7) at (-7,0) {};
		
		\foreach \i in {1,2,3} {
			\node[vtx] (r\i) at (\i,0) {};
		}
		\node[wvtx] (r4) at (4,0) {};
		
		\node[vtx] (u1) at (0,1) {};
		
		\draw (l7)--(l6)--(l5)--(l4)--(l3)--(l2)--(l1)--(c)
		--(r1)--(r2)--(r3)--(r4);
		\draw (c)--(u1);
		
		\node[
		draw,
		dashed,
		rounded corners,
		inner sep=6pt,
		fit=(l1)(r1)(u1)(c)
		] {};
		
		\node[
		draw,
		dashed,
		rounded corners,
		inner sep=6pt,
		fit=(l2)(l4)
		] {};
		
		\node[
		draw,
		dashed,
		rounded corners,
		inner sep=6pt,
		fit=(l5)(l7)
		] {};
		
		\node[
		draw,
		dashed,
		rounded corners,
		inner sep=6pt,
		fit=(r2)(r4)
		] {};
		
	\end{tikzpicture}
	
	\caption{The graph $G_{2,4}=S(7,4,1)$ with
		$\tr(G_{2,4})=24=6\cdot\rho_{2}(G_{2,4})$. The dashed regions are the
		central part and the three additional blocks, each of diameter at most
		$2$. The left leg contains two such blocks, corresponding to $q_1=2$.
		The white vertices form a $2$-packing of cardinality $4$.}
	\label{fig-Gkh}
\end{figure}
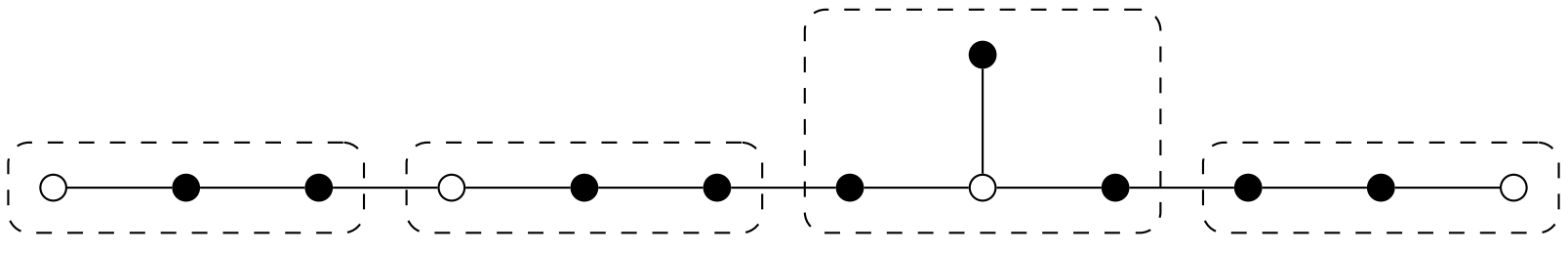

Figure~\ref{fig-Gkh} contains the graph $G_{k,h}$ from the proof of
Theorem~\ref{thm-packing} for $k=2$, $h=4$, and
$(q_1,q_2,q_3)=(2,1,0)$.

\section{Minimum degree and domination}\label{sec-5}

The packing bound from Theorem~\ref{thm-packing} immediately yields
estimates in terms of several classical graph parameters. For $k=1$,
a $1$-packing is simply an independent set, and hence
$\rho_1(G)=\alpha(G)$. For $k=2$, standard bounds on the $2$-packing
number in terms of the minimum degree and the domination number give
the corresponding estimates for the triameter. We therefore obtain
the following corollary.

\begin{corollary}\label{cor-degree}
For any connected graph $G$ with $n\geq 3$ vertices, we have
\begin{align}
\tr(G)&\leq 4\alpha(G)-1, \label{eq-alpha}\\
\tr(G)&\leq 6\rho_{2}(G)\leq 6\floor*{\frac{n}{\delta(G)+1}}, \label{eq-delta}\\
\tr(G)&\leq 6\gamma(G). \label{eq-gamma}
\end{align}
Moreover, the bounds in terms of $\alpha(G)$, $\rho_{2}(G)$ and $\gamma(G)$ are tight for every positive integer value of the corresponding parameter. In particular, Conjecture~\ref{conj-DA} holds.
\end{corollary}
\begin{proof}
The inequality~\eqref{eq-alpha} is Theorem~\ref{thm-packing} for $k=1$, since $\rho_{1}(G)=\alpha(G)$, and the first inequality in~\eqref{eq-delta} is Theorem~\ref{thm-packing} for $k=2$.

Let $P$ be a $2$-packing in $G$. Then the closed neighborhoods $N_{G}[x]$, $x\in P$, are pairwise disjoint, and each of them contains at least $\delta(G)+1$ vertices. Thus, $|P|(\delta(G)+1)\leq n$, which implies $\rho_{2}(G)\leq\floor*{\frac{n}{\delta(G)+1}}$.

Let $D$ be a dominating set of $G$. For every $x\in P$, choose a vertex of $D\cap N_G[x]$. Since the closed neighborhoods $N_G[x]$, $x\in P$, are pairwise disjoint, these chosen vertices are distinct. Hence, $|P|\leq |D|$, and therefore $\rho_2(G)\leq\gamma(G)$. Together with the first inequality in~\eqref{eq-delta}, this proves~\eqref{eq-gamma}.

To see that the bound in terms of $\alpha(G)$ is sharp, consider the
graphs $G_{1,h}$ from the proof of Theorem~\ref{thm-packing}. Since
$\rho_{1}(G)=\alpha(G)$, these graphs satisfy
\[
\alpha(G_{1,h})=h
\qquad\text{and}\qquad
\tr(G_{1,h})=4h-1.
\]
For example, choosing $(q_{1},q_{2},q_{3})=(h-1,0,0)$ yields
$U(2h-2,0,0)$, that is, a triangle with a pendant path of length
$2h-2$.

For the bounds in terms of $\rho_{2}(G)$ and $\gamma(G)$, consider the
graphs $G_{2,h}$ from the proof of Theorem~\ref{thm-packing}. These are
spiders whose $i$-th leg has length $3q_i+1$, and
\[
\tr(G_{2,h})=6h=6\rho_{2}(G_{2,h}).
\]
We claim that $\gamma(G_{2,h})=h$. Indeed, the center together with the
vertices at depths $3,6,\ldots,3q_i$ on the $i$-th leg, for
$i=1,2,3$, forms a dominating set of cardinality $h$. The center
dominates all vertices at depth $1$, while a vertex at depth $3j$
dominates the vertices at depths $3j-1$, $3j$, and $3j+1$ on its leg.
Thus
\[
\rho_{2}(G_{2,h})\leq \gamma(G_{2,h})\leq h.
\]
Since $\rho_{2}(G_{2,h})=h$, it follows that
\[
\gamma(G_{2,h})=h,
\]
and hence
\[
\tr(G_{2,h})
=6\rho_{2}(G_{2,h})
=6\gamma(G_{2,h}).
\]
Thus both bounds are sharp.
\end{proof}

The inequality~\eqref{eq-gamma} recovers Proposition~\ref{prop-Das-dom}. The next proposition shows that the coefficient $6$ in~\eqref{eq-delta} cannot be improved.

\begin{proposition}\label{prop-six}
	For any fixed integer $\delta\geq 2$, the coefficient $6$ in the bound
	\[
	\tr(G)\leq\frac{6n}{\delta(G)+1}
	\]
	is asymptotically best possible.
\end{proposition}
\begin{proof}
	Fix $j\in\mathbb{N}$ and consider the path with vertices
	$0,1,\dots,3j$. Construct its clique blow-up $G$ as follows. The bags
	corresponding to $0$ and $3j$ have $\delta$ vertices, the bags
	corresponding to $3,6,\dots,3j-3$ have $\delta-1$ vertices, and all
	remaining bags consist of a single vertex. Thus,
	\[
	|V(G)|=2\delta+(j-1)(\delta-1)+2j=(j+1)(\delta+1).
	\]
	
	A vertex in either of the two end bags has degree
	$(\delta-1)+1=\delta$. A vertex in the bag corresponding to $3i$,
	where $1\leq i\leq j-1$, has degree $(\delta-2)+2=\delta$, while every
	remaining vertex has degree at least $\delta$. Hence $\delta(G)=\delta$.
	
	For vertices belonging to distinct bags, their distance in $G$ equals
	the distance between the corresponding vertices of the underlying
	path. Therefore, a triple whose vertices belong to three distinct bags
	has perimeter at most $6j$, and a triple contained in a single bag
	has perimeter $3$.
	
	Now suppose that exactly two vertices of the triple belong to the same bag.
	Their mutual distance is $1$, and if this bag corresponds to the
	vertex $r\in\{0,\dots,3j\}$, while the third vertex belongs to the bag
	corresponding to $s$, then the perimeter is
	\[
	1+2|r-s|\leq 1+6j.
	\]
	Equality is possible only when $\{r,s\}=\{0,3j\}$. Since both end bags
	contain at least two vertices, we may take two vertices from the bag
	of $0$ and one vertex from the bag of $3j$. This triple has perimeter
	$6j+1$. Consequently,
	\[
	\tr(G)=6j+1.
	\]
	
	It follows that
	\[
	\frac{\tr(G)(\delta(G)+1)}{|V(G)|}
	=\frac{6j+1}{j+1}\longrightarrow 6
	\qquad\text{as }j\to\infty.
	\qedhere
	\]
\end{proof}

The estimate in~\eqref{eq-delta} is obtained by bounding the size of a
$2$-packing using pairwise disjoint closed neighborhoods. This argument
uses only the first neighborhood of each packing vertex and therefore
does not take into account any further local expansion of the graph.
For graph classes in which metric balls grow sufficiently fast, one can
apply the same idea to $2s$-packings and balls of radius $s$, obtaining
stronger bounds. To make this precise, for $s\in\mathbb{N}$ define
\[
b_s(G)=\min_{v\in V(G)}|B_G(v,s)|.
\]
The following proposition gives the corresponding refinement.

\begin{proposition}\label{prop-ball-bound}
	Let $G$ be a connected graph of order $n$, and let $s\in\mathbb{N}$. Then
	\[
	\tr(G)\leq
	(4s+2)\left\lfloor\frac{n}{b_s(G)}\right\rfloor+2s-2.
	\]
\end{proposition}
\begin{proof}
	Let $P$ be a maximum $2s$-packing in $G$. The balls
	$B_G(v,s)$, $v\in P$, are pairwise disjoint, and hence
	\[
	|P|\,b_s(G)\leq n.
	\]
	Therefore,
	\[
	\rho_{2s}(G)\leq
	\left\lfloor\frac{n}{b_s(G)}\right\rfloor.
	\]
	Applying Theorem~\ref{thm-packing} with $k=2s$ gives
	\[
	\tr(G)\leq
	(4s+2)\rho_{2s}(G)+2s-2
	\leq
	(4s+2)\left\lfloor\frac{n}{b_s(G)}\right\rfloor+2s-2.
	\]
\end{proof}

For graphs without short cycles, the quantity $b_2(G)$ can be bounded
from below in terms of the minimum degree, which yields the following
explicit estimate.

\begin{corollary}\label{cor-girth-five}
	Let $G$ be a connected graph of order $n$ and girth at least five. Then
	\[
	\tr(G)\leq
	10\left\lfloor\frac{n}{\delta(G)^2+1}\right\rfloor+2.
	\]
\end{corollary}
\begin{proof}
	Let $v\in V(G)$. Since $G$ has girth at least five, no two neighbors
	of $v$ are adjacent, and the neighborhoods of distinct neighbors of $v$
	are disjoint outside $v$. Each neighbor of $v$ has at least
	$\delta(G)-1$ neighbors different from $v$, and hence
	\[
	|B_G(v,2)|
	\geq
	1+\deg(v)+\deg(v)(\delta(G)-1)
	\geq
	\delta(G)^2+1.
	\]
	Thus $b_2(G)\geq \delta(G)^2+1$. Applying
	Proposition~\ref{prop-ball-bound} with $s=2$, we obtain
	\[
	\tr(G)\leq
	10\left\lfloor\frac{n}{b_2(G)}\right\rfloor+2
	\leq
	10\left\lfloor\frac{n}{\delta(G)^2+1}\right\rfloor+2.
	\]
\end{proof}

\section{Vertex connectivity}\label{sec-6}

In this section, we study how large triameter can be in a graph of
prescribed vertex connectivity. Intuitively, a large triameter requires
three vertices to be separated by long metric branches, whereas high
vertex connectivity prevents these branches from being joined through
small vertex cuts. The canonical projection introduced in
Section~\ref{sec-3} makes this tension explicit: suitable separators in
the corresponding spider or triangle-with-legs model pull back to
vertex cuts of the original graph.

Using this idea, we derive a sharp upper bound on the triameter in terms
of the order and vertex connectivity, prove
Theorem~\ref{thm-connectivity}, and determine the extremal graphs. We
keep the notation of Section~\ref{sec-3}. The key step is the following
lemma.

\begin{lemma}\label{lem-separators}
	Let $G$ be a connected graph, let
	$T=\{a_1,a_2,a_3\}\subseteq V(G)$, and put
	\[
	p=d_G(a_1,a_2,a_3).
	\]
	Then $G$ contains
	\[
	\floor*{\frac{p-3}{2}}
	\]
	pairwise disjoint $T$-separators. Moreover, these separators can be
	chosen as fibers $\varphi^{-1}(x)$ of vertices $x\in V(H_T)$.
\end{lemma}
\begin{proof}
	By Lemma~\ref{lem-transfer}, it is enough to find
	$\floor*{(p-3)/2}$ vertices of $H_T$, distinct from the feet
	$A_1,A_2,A_3$, whose deletion separates the triple of feet.
	The corresponding fibers under $\varphi$ are then $T$-separators in
	$G$, and fibers of distinct vertices are pairwise disjoint.
	
	We distinguish the two possible types of the canonical model $H_T$.
	
	Suppose first that $p$ is even. Then $H_T$ is a tree and
	\[
	|V(H_T)|=\frac p2+1.
	\]
	If $H_T$ has three leaves, every non-leaf vertex separates at least one
	foot from the other two. Hence all
	\[
	|V(H_T)|-3=\frac p2-2
	=\floor*{\frac{p-3}{2}}
	\]
	non-leaf vertices may be chosen.
	
	If $H_T$ is a path, one of the feet lies in its interior. Every interior
	vertex other than this foot separates the three feet into at least two
	components. There are again
	\[
	\frac p2-2
	=\floor*{\frac{p-3}{2}}
	\]
	such vertices.
	
	Now suppose that $p$ is odd. Then
	$H_T=U(\ell_1,\ell_2,\ell_3)$. On the $i$th leg, every vertex of depth
	$0,1,\ldots,\ell_i-1$ separates the foot $A_i$ from the other two feet.
	Thus we may choose
	\[
	\ell_1+\ell_2+\ell_3
	=\frac{p-3}{2}
	\]
	vertices in total.
	
	In every case, this gives
	$\floor*{(p-3)/2}$ vertices of $H_T$ whose fibers are pairwise disjoint
	$T$-separators in $G$.
\end{proof}

\begin{remark}
	The bound $\floor*{\frac{p-3}{2}}$ in
	Lemma~\ref{lem-separators} is best possible. To see why no larger
	universal bound can hold, fix a geodesic $P_{ij}$ between $a_i$ and
	$a_j$ for each $1\leq i<j\leq 3$. A $T$-separator must separate at
	least two of the three pairs $(a_1,a_2)$, $(a_1,a_3)$, $(a_2,a_3)$.
	Consequently, it must contain an internal vertex of at least two of the
	three geodesics $P_{12},P_{13},P_{23}$.
	
	The three geodesics have
	\[
	d(a_1,a_2)+d(a_1,a_3)+d(a_2,a_3)-3=p-3
	\]
	internal vertices in total, counted with multiplicity. Hence, if
	$S_1,\ldots,S_r$ are pairwise disjoint $T$-separators, each $S_t$
	accounts for at least two of these occurrences, and therefore
	\[
	2r\leq p-3.
	\]
	Thus
	\[
	r\leq \floor*{\frac{p-3}{2}},
	\]
	showing that the number of pairwise disjoint $T$-separators provided by
	Lemma~\ref{lem-separators} is optimal.
	
	There is a parallel edge version. Using the edges of the canonical
	model $H_T$, one obtains $\floor*{\frac{p}{2}}$ pairwise disjoint edge
	sets, each of which separates at least two of the three pairs
	$a_1,a_2,a_3$; the same counting argument shows that this number is
	optimal. When only two distinguished vertices are considered, the
	corresponding statement reduces to the classical fact that the distance
	between two vertices equals the maximum number of pairwise disjoint
	edge cuts separating them~\cite{Robacker:56,Schrijver:03}.
\end{remark}

We are now ready to prove the main result of this section. The idea is
to apply Lemma~\ref{lem-separators} to a triametral triple and compare
the resulting family of pairwise disjoint separators with the vertex
connectivity of the graph.

\begin{proof}[Proof of Theorem~\ref{thm-connectivity}]
\textbf{Upper bound.} Let $G$ be a connected graph with $n$ vertices, let $a_{1},a_{2},a_{3}$ be its triametral triple, and put $T=\{a_1,a_2,a_3\}$. By Lemma~\ref{lem-separators}, $G$ contains $\floor*{\frac{\tr(G)-3}{2}}$ pairwise disjoint $T$-separators. Each of them is a vertex cut, so it contains at least $\kappa(G)$ vertices, and none of them contains $a_{1},a_{2},a_{3}$. Therefore,
\begin{align}\label{eq-kappa}
\kappa(G)\floor*{\frac{\tr(G)-3}{2}}\leq n-3.
\end{align}
Put $q=\floor*{\frac{n-3}{\kappa(G)}}$. Then~\eqref{eq-kappa} gives $\floor*{\frac{\tr(G)-3}{2}}\leq q$, which means $\frac{\tr(G)-3}{2}<q+1$. Hence, $\tr(G)\leq 2q+4$, that is,
\begin{align}\label{eq-kappa-2}
\tr(G)\leq 2\floor*{\frac{n-3}{\kappa(G)}}+4.
\end{align}

\textbf{Constructions.} If $n=k+1$, then $G=K_{n}$ is the only graph with $\kappa(G)=k$, and $\tr(K_{n})=3$. If $n=k+2$, then the graph $K_{n}-e$ obtained by deleting one edge has connectivity $k$ and triameter $4$, which agrees with~\eqref{eq-kappa-2}.

Now let $n\geq k+3$. Put $q=\floor*{\frac{n-3}{k}}\geq 1$ and $r=n-(kq+3)$, so that $0\leq r<k$. Consider the spider $S(q,1,1)$. It has $q$ non-leaf vertices: the center and $q-1$ interior vertices of the long leg. Construct its clique blow-up $G_{n,k}$ in which the bag of every non-leaf vertex has $k$ vertices, the bag of the foot of the long leg and the bag of one of the short feet consist of a single vertex, and the bag of the other short foot has $r+1$ vertices (see Figure~\ref{fig-Gnk}). Then $|V(G_{n,k})|=kq+3+r=n$.

\begin{figure}[ht]
\centering
\begin{tikzpicture}[scale=1,every node/.style={font=\small}]
\node[vtx,label=below:$x$] (x) at (-4.5,0) {};
\node[draw,ellipse,minimum width=1.3cm,minimum height=0.8cm] (C1) at (-2.4,0) {$K_{k}$};
\node[draw,ellipse,minimum width=1.3cm,minimum height=0.8cm] (C0) at (0,0) {$K_{k}$};
\node[vtx,label=right:$y$] (y) at (2.2,1) {};
\node[draw,ellipse,minimum width=1.3cm,minimum height=0.8cm] (Z) at (2.4,-1) {$K_{r+1}$};
\draw[line width=2pt] (x)--(C1)--(C0)--(y);
\draw[line width=2pt] (C0)--(Z);
\end{tikzpicture}
\caption{The graph $G_{n,k}$ for $q=2$. Thick lines denote complete joins between the corresponding bags.}
\label{fig-Gnk}
\end{figure}
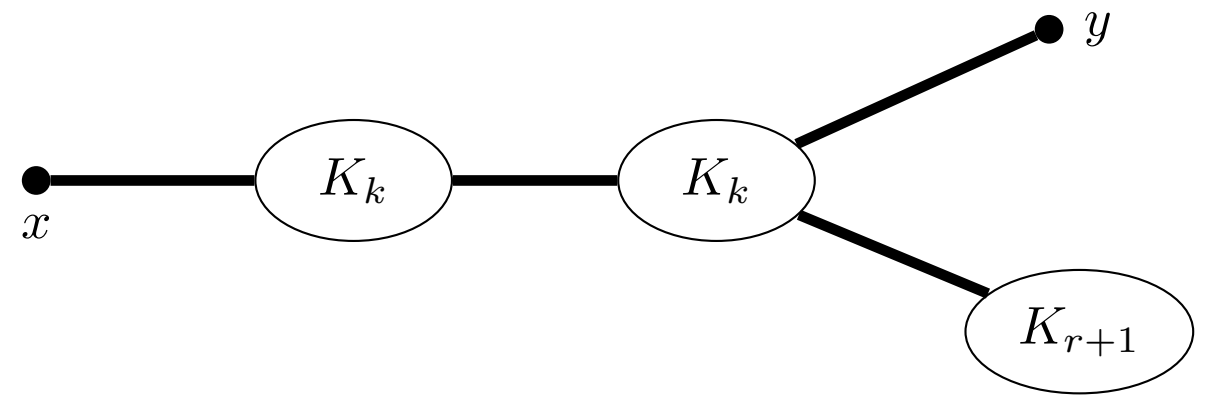

\textbf{Claim 1.} $\kappa(G_{n,k})=k$.

Let $X\subseteq V(G_{n,k})$ with $|X|<k$. Every bag of a non-leaf vertex has $k$ vertices, so it contains a vertex outside $X$. The remaining vertices of the bags of the non-leaf vertices induce a connected subgraph, since the bags of adjacent vertices are completely joined. Every remaining vertex of a leaf bag is adjacent to all vertices of the neighboring non-leaf bag. Therefore, $G_{n,k}-X$ is connected. On the other hand, deleting the bag of the center (which has $k$ vertices) disconnects the graph into three parts corresponding to the three legs. This proves the claim.

\textbf{Claim 2.} $\tr(G_{n,k})=2q+4$.

The distance between two vertices of $G_{n,k}$ from different bags equals the distance between the corresponding vertices of $S(q,1,1)$. Therefore, a triple of vertices from three distinct bags has perimeter at most $\tr(S(q,1,1))=2q+4$, and the three feet attain this value. A triple using some bag twice has perimeter at most $1+2(q+1)=2q+3$, and a triple from a single bag has perimeter $3$. This proves the claim.

By Claims~1 and~2, the bound~\eqref{eq-kappa-2} is attained for all $n\geq k+3$, which finishes the proof.
\end{proof}

Since $\frac{3(n-2)}{k}+3-\left(\frac{2(n-3)}{k}+4\right)=\frac{n}{k}-1>0$, Theorem~\ref{thm-connectivity} indeed strengthens Proposition~\ref{prop-Das-kappa}. Also, for $k=1$ the graph $G_{n,1}$ is a tree with three leaves, in accordance with Proposition~\ref{prop-Das-order}.

\subsection{Extremal graphs}

We now describe the $k$-connected graphs attaining the bound~\eqref{eq-kappa-2}.

\begin{proposition}\label{prop-equality}
Let $k\geq 2$, $n\geq k+2$ and $s=\floor*{\frac{n-3}{k}}+2$. A $k$-connected graph $G$ with $n$ vertices has $\tr(G)=2s$ if and only if $G$ is a spanning subgraph of a clique blow-up of one of the following graphs:
\begin{enumerate}
\item a tree with exactly three leaves and $s$ edges, in which the bags of all non-leaf vertices have at least $k$ vertices;
\item a path with $s$ edges, in which the bags of all interior vertices have exactly $k$ vertices and the bags of both ends consist of a single vertex; this case occurs only if $n=k(s-1)+2$.
\end{enumerate}
\end{proposition}
\begin{proof}
\textbf{Necessity.} Let $G$ be a $k$-connected graph with $\tr(G)=2s$. Fix a triametral triple $a_{1},a_{2},a_{3}$ and consider the projection $\varphi\colon V(G)\to V(H_{T})$. Since $p=2s$ is even, $H_{T}$ is a tree with three leaves or a path, and it has $s$ edges by Lemma~\ref{lem-projection}.

Let $x\in V(H_{T})$ be a vertex whose deletion disconnects $H_{T}$. Then the feet do not all lie in one connected component of $H_{T}-x$ (this includes the foot lying inside the path, whose deletion separates the other two feet). By the argument from the second statement of Lemma~\ref{lem-transfer}, $\varphi^{-1}(x)$ is a vertex cut in $G$, and hence $|\varphi^{-1}(x)|\geq k$.

Since $\varphi$ is non-expanding, every edge of $G$ joins two vertices from one preimage or from the preimages of two adjacent vertices of $H_{T}$. Thus, $G$ is a spanning subgraph of the clique blow-up of $H_{T}$ whose bags are the preimages $\varphi^{-1}(x)$, $x\in V(H_{T})$. If $H_{T}$ is a tree with three leaves, then all its non-leaf vertices disconnect it, and we obtain the first case. If $H_{T}$ is a path, then all its $s-1$ interior vertices disconnect it, and hence $n\geq k(s-1)+2$. On the other hand, the definition of $s$ gives $n-3<k(s-1)$, that is, $n\leq k(s-1)+2$. Therefore, $n=k(s-1)+2$, and all the bag sizes are uniquely determined, which is the second case.

\textbf{Sufficiency.} Let $\widehat{F}$ be a clique blow-up of a tree $F$ as in the first case. The distance between vertices from different bags of $\widehat{F}$ equals the distance between the corresponding vertices of $F$. Hence, three vertices from the three leaf bags give perimeter $\tr(F)=2s$, and no triple from three distinct bags gives more. The diameter of a tree with three leaves and $s$ edges is at most $s-1$, so a triple using some bag twice has perimeter at most $1+2(s-1)$. Thus, $\tr(\widehat{F})=2s$. Similarly, if $\widehat{F}$ is a clique blow-up of a path as in the second case, then the two end vertices together with any third vertex give perimeter $2s$. A repeated bag must be interior, and its distance to any other bag is at most $s-1$. So again $\tr(\widehat{F})=2s$.

Now let $G$ be a $k$-connected spanning subgraph of $\widehat{F}$. Deleting edges does not decrease distances, so $\tr(G)\geq\tr(\widehat{F})=2s$. On the other hand, $\tr(G)\leq 2s$ by~\eqref{eq-kappa-2}. Therefore, $\tr(G)=2s$.
\end{proof}

For the smallest order $n=k(s-2)+3$, the structure of the extremal graphs is much more rigid.

\begin{theorem}\label{thm-Kk3}
Let $k\geq 2$, $s\geq 3$ and $n=k(s-2)+3$. A $k$-connected graph $G$ with $n$ vertices has $\tr(G)=2s$ if and only if $G$ is a spanning subgraph of a clique blow-up of a tree with three leaves and $s$ edges, in which the bags of all non-leaf vertices have exactly $k$ vertices and the bags of the leaves consist of a single vertex. Moreover, for any triametral triple $a_{1},a_{2},a_{3}$ of such a graph $G$:
\begin{enumerate}
\item there is such a blow-up in which $\{a_{1}\}$, $\{a_{2}\}$, $\{a_{3}\}$ are the bags of the leaves;
\item $G$ contains a spanning subdivision of $K_{k,3}$ in which $a_{1},a_{2},a_{3}$ form the part of size three, and all the $k$ paths starting at $a_{i}$ have the same length $\ell_{i}\geq 1$, where $\ell_{1}+\ell_{2}+\ell_{3}=s$.
\end{enumerate}
\end{theorem}
\begin{proof}
\textbf{Necessity.} Note that here $\floor*{\frac{n-3}{k}}+2=s$. Let $\tr(G)=2s$, fix a triametral triple $a_{1},a_{2},a_{3}$, and consider the projection $\varphi\colon V(G)\to V(H_{T})$. If $H_{T}$ were a path, then by Proposition~\ref{prop-equality} we would have $n=k(s-1)+2>k(s-2)+3$, a contradiction. Thus, $H_{T}=S(\ell_{1},\ell_{2},\ell_{3})$ with $\ell_{1},\ell_{2},\ell_{3}\geq 1$ and $\ell_{1}+\ell_{2}+\ell_{3}=s$. By the proof of Proposition~\ref{prop-equality}, the preimages of the $s-2$ non-leaf vertices of $H_{T}$ have at least $k$ vertices each, and the preimages of the three feet contain $a_{1},a_{2},a_{3}$, respectively. Since $n=k(s-2)+3$, all these inequalities are equalities: the preimage of every non-leaf vertex has exactly $k$ vertices, and $\varphi^{-1}(A_{i})=\{a_{i}\}$. This proves the necessity together with the first statement.

\textbf{Sufficiency} follows from Proposition~\ref{prop-equality}.

It remains to prove the second statement. We use the blow-up from the first statement, whose bags are the preimages $\varphi^{-1}(x)$, $x\in V(H_{T})$.

\textbf{Claim 1.} Let $X,Y$ be the bags of two adjacent non-leaf vertices on the $i$--th leg of $H_{T}$, where $X$ is closer to $a_{i}$ (here $Y$ may be the bag of the center). Then the edges of $G$ between $X$ and $Y$ contain a perfect matching.

To the contrary, suppose that there is no such perfect matching. Since $|X|=|Y|=k$, Hall's theorem~\cite{Hall:35} gives a set $X'\subseteq X$ with $|N_{G}(X')\cap Y|<|X'|$. Put $Z=(X\setminus X')\cup(N_{G}(X')\cap Y)$. Then $|Z|=k-|X'|+|N_{G}(X')\cap Y|<k$. Since the edges of $G$ join only vertices from equal or adjacent bags, in $G-Z$ there are no edges between the part of the $i$--th leg closer to $a_{i}$ than $Y$ and the set $Y\setminus Z$. The vertex $a_{i}$ survives in $G-Z$, and so does some vertex of $Y$, since $|N_{G}(X')\cap Y|<k$. Therefore, $G-Z$ is disconnected, which contradicts $\kappa(G)\geq k$. The obtained contradiction proves the claim.

\textbf{Claim 2.} The vertex $a_{i}$ is adjacent to all vertices of the bag of its neighbor in $H_{T}$.

Indeed, all neighbors of $a_{i}$ lie in this bag, which has $k$ vertices, and $\deg_{G}(a_{i})\geq\kappa(G)\geq k$.

By Claims~1 and~2, starting at $a_{i}$ and following the perfect matchings along the $i$--th leg, we obtain $k$ internally disjoint paths of length $\ell_{i}$ from $a_{i}$ to the $k$ distinct vertices of the bag of the center. These paths cover all the vertices of the $i$--th leg. The union of these paths over $i=1,2,3$ is the required spanning subdivision of $K_{k,3}$, in which the bag of the center is the part of size $k$.
\end{proof}

\begin{remark}
For $k\geq 3$, the spanning subdivision of $K_{k,3}$ from Theorem~\ref{thm-Kk3} is not $k$-connected itself, unless $k=3$ and $\ell_{1}=\ell_{2}=\ell_{3}=1$ (in which case it is $K_{3,3}$). Indeed, its vertices of the part of size $k$ have degree $3$, and its interior path vertices have degree $2$. Thus, the other edges of $G$ are essential for the $k$-connectivity.
\end{remark}

For $k=2$, Theorem~\ref{thm-Kk3} gives the following uniqueness result.

\begin{corollary}\label{cor-biconnected}
For any biconnected graph $G$ with $n$ vertices, we have
\[\tr(G)\leq 2\floor*{\frac{n+1}{2}},\]
and this bound is tight for all $n\geq 4$. Moreover, if $\tr(G)=n+1$, then $n$ is odd, $G$ has a unique triametral triple, and $G$ contains a spanning theta subgraph whose three paths have even lengths; the unique triametral triple consists of the midpoints of these three paths.
\end{corollary}
\begin{proof}
The bound is~\eqref{eq-kappa-2} for $\kappa(G)\geq 2$, since $\floor*{\frac{n-3}{2}}+2=\floor*{\frac{n+1}{2}}$. For even $n$, it is attained by the cycle $C_{n}$. For odd $n\geq 5$, take a theta graph with paths of lengths $2\ell_{1},2\ell_{2},2\ell_{3}$, where $\ell_{1},\ell_{2},\ell_{3}\geq 1$ and $\ell_{1}+\ell_{2}+\ell_{3}=\frac{n+1}{2}$. It has $2+(2\ell_{1}-1)+(2\ell_{2}-1)+(2\ell_{3}-1)=n$ vertices, and the midpoints of its paths are pairwise at distances $\ell_{i}+\ell_{j}$, which gives perimeter $n+1$ (see Figure~\ref{fig-theta}).

\begin{figure}[ht]
\centering
\begin{tikzpicture}[scale=1,every node/.style={font=\small}]
\node[wvtx,label=above:$h_1$] (h1) at (0,1.6) {};
\node[wvtx,label=below:$h_2$] (h2) at (0,-1.6) {};
\node[vtx,label=left:$a_1$] (m1) at (-1.6,0) {};
\node[vtx,label=left:$a_2$] (m2) at (0,0) {};
\node[wvtx] (p1) at (1.3,0.9) {};
\node[vtx,label=right:$a_3$] (m3) at (1.9,0) {};
\node[wvtx] (p2) at (1.3,-0.9) {};
\draw (h1)--(m1)--(h2);
\draw (h1)--(m2)--(h2);
\draw (h1)--(p1)--(m3)--(p2)--(h2);
\end{tikzpicture}
\caption{A theta graph with $n=7$ vertices and paths of lengths $2,2,4$. The midpoints $a_1,a_2,a_3$ form its unique triametral triple, with perimeter $2+3+3=8=n+1$.}
\label{fig-theta}
\end{figure}
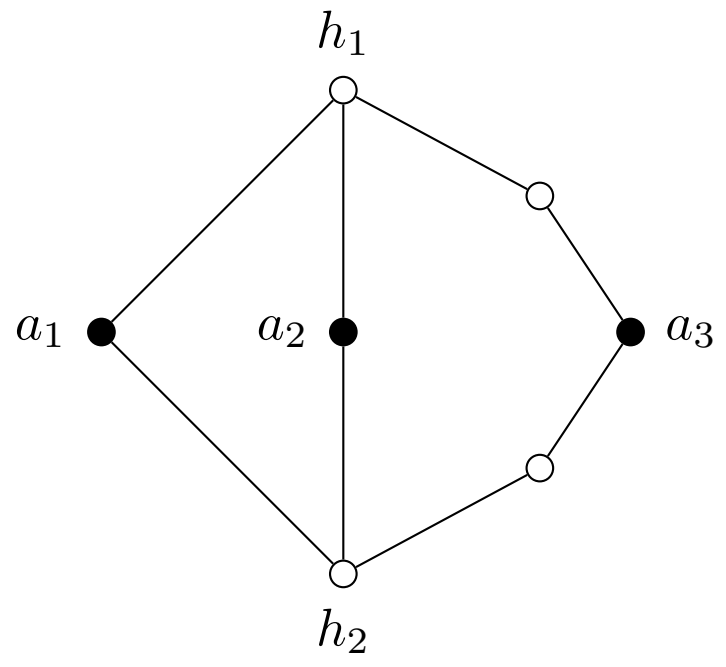

Now suppose that $\tr(G)=n+1=2s$. Then $n$ is odd, $n\geq 5$ (since $K_{3}$ is the only biconnected graph with three vertices), $s\geq 3$ and $n=2(s-2)+3$. Fix a triametral triple $a_{1},a_{2},a_{3}$. By Theorem~\ref{thm-Kk3} for $k=2$, the graph $G$ contains a spanning theta subgraph $\Theta$ with hubs $h_{1},h_{2}$, whose paths have lengths $2\ell_{1},2\ell_{2},2\ell_{3}$ and contain $a_{1},a_{2},a_{3}$ as their midpoints.

We claim that the midpoints form the unique triple with perimeter $2s$ in $\Theta$. Consider a triple of vertices $v_{1},v_{2},v_{3}$ in $\Theta$. If it meets at most two of the three paths of $\Theta$, then it lies on a cycle of length $2(\ell_{i}+\ell_{j})<2s$, and its perimeter in $\Theta$ is at most the length of this cycle. This includes the triples containing a hub. Otherwise, we can assume that $v_{i}$ is an interior vertex of the $i$--th path, and write $d_{\Theta}(h_{1},v_{i})=\ell_{i}+u_{i}$ with $|u_{i}|<\ell_{i}$. Considering the routes through $h_{1}$ and through $h_{2}$, we obtain
\[d_{\Theta}(v_{i},v_{j})\leq\min\{\ell_{i}+u_{i}+\ell_{j}+u_{j},\ \ell_{i}-u_{i}+\ell_{j}-u_{j}\}=\ell_{i}+\ell_{j}-|u_{i}+u_{j}|.\]
Summing over the three pairs gives
\[d_{\Theta}(v_{1},v_{2},v_{3})\leq 2s-|u_{1}+u_{2}|-|u_{1}+u_{3}|-|u_{2}+u_{3}|.\]
The equality $d_{\Theta}(v_{1},v_{2},v_{3})=2s$ forces $u_{1}+u_{2}=u_{1}+u_{3}=u_{2}+u_{3}=0$, which implies $u_{1}=u_{2}=u_{3}=0$. This proves our claim.

Finally, since $\Theta$ is a spanning subgraph of $G$, we have $d_{G}(u,v)\leq d_{\Theta}(u,v)$ for all $u,v\in V(G)$. Therefore, no triple other than $a_{1},a_{2},a_{3}$ has perimeter $2s$ in $G$.
\end{proof}

Note that the theta subgraph from Corollary~\ref{cor-biconnected} need not be induced or unique. Also, for even $n$, the equality $\tr(G)=n$ does not force a unique triametral triple, as even cycles show.

\section{Open questions}\label{sec-7}

We conclude the paper with several open problems suggested by our results.

\textbf{Question 1.} Theorem~\ref{thm-packing} determines the distances between the vertices of a triametral triple in graphs attaining the equality in~\eqref{eq-packing}, but not the whole graph. Characterize the graphs attaining this equality, at least for $k=1$, that is, the graphs $G$ with $\tr(G)=4\alpha(G)-1$.

\textbf{Question 2.} Find the exact value of $\max\{\tr(G):|V(G)|=n,\ \delta(G)=\delta\}$. By Proposition~\ref{prop-six}, the coefficient $6$ in~\eqref{eq-delta} is asymptotically correct, but the exact maximum generally differs from $\frac{6n}{\delta+1}$.

\textbf{Question 3.} Do Theorem~\ref{thm-packing} and Lemma~\ref{lem-separators} have analogs for $r\geq 4$ prescribed vertices, with the sum of all $\binom{r}{2}$ distances in place of the perimeter? Note that for $r\geq 4$ the tight span of an $r$-point metric space is no longer a tree, so a different projection is needed.


\section*{Declaration of generative AI use}
An initial version of the results and proofs presented here was produced with the assistance of the AI system GPT-Astra.
The AI system Claude Opus 5.5 was subsequently used to check the arguments step by step, to suggest corrections and additions, and to help restructure the material and prepare the manuscript.
All statements and proofs in the final version were independently verified by the human authors, who selected, revised and rewrote the material as needed.
The authors take full responsibility for the content of this article.



\begin{thebibliography}{99}

\bibitem{BanChe:08}
H.-J.~Bandelt, V.~Chepoi, Metric graph theory and geometry: a survey, in: \textit{Surveys on Discrete and Computational Geometry}, Contemp. Math. \textbf{453}, Amer. Math. Soc., Providence, RI, 2008, pp.~49--86.

\bibitem{Das:21}
A.~Das, Triameter of graphs, \textit{Discuss. Math. Graph Theory} \textbf{41} (2021) 601--616.

\bibitem{DasAbd:25}
A.~Das, C.~Abdio\u{g}lu, The triameter of a graph, \textit{Necmettin Erbakan Univ. J. Sci. Eng.} \textbf{7(3)} (2025) 535--543.

\bibitem{Dress:84}
A.~W.~M.~Dress, Trees, tight extensions of metric spaces, and the cohomological dimension of certain groups: a note on combinatorial properties of metric spaces, \textit{Adv. Math.} \textbf{53} (1984) 321--402.

\bibitem{HakKozOl:22}
A.~Hak, S.~Kozerenko, B.~Oliynyk, A note on the triameter of graphs, \textit{Discrete Appl. Math.} \textbf{309} (2022) 278--284.

\bibitem{Hall:35}
P.~Hall, On representatives of subsets, \textit{J. London Math. Soc.} \textbf{10} (1935) 26--30.

\bibitem{HenYeo:14}
M.~A.~Henning, A.~Yeo, A new lower bound for the total domination number in graphs
proving a Graffiti.pc conjecture,
\textit{Discrete Appl. Math.} \textbf{173} (2014) 45--52.

\bibitem{Isbell:64}
J.~R.~Isbell, Six theorems about injective metric spaces, \textit{Comment. Math. Helv.} \textbf{39} (1964) 65--76.

\bibitem{KolaPan:15}
S.~R.~Kola, P.~Panigrahi, A lower bound for radio $k$-chromatic number of an arbitrary graph,
\textit{Contrib. Discrete Math.} \textbf{10} (2015) 45--56.

\bibitem{Koz:26}
S.~Kozerenko, Graphs with small triameters, \textit{Discrete Math. Lett.} \textbf{18} (2026) 55--62.

\bibitem{Laurent:96}
M.~Laurent, Graphic vertices of the metric polytope, \textit{Discrete Math.} \textbf{151} (1996) 131--153.

\bibitem{Robacker:56}
J.~T.~Robacker, Min-max theorems on shortest chains and disjoint cuts of a network, Research Memorandum RM-1660, The RAND Corporation, Santa Monica, 1956.

\bibitem{SahaBas:20}
L.~Saha, A.~R.~Basunia,
Radio graceful labelling of graphs,
\textit{Theory Appl. Graphs} \textbf{7} (2020) Article~7.

\bibitem{SahaPan:15}
L.~Saha, P.~Panigrahi,
A lower bound for radio $k$-chromatic number,
\textit{Discrete Appl. Math. }\textbf{192} (2015) 87--100.

\bibitem{Schrijver:03}
A.~Schrijver, \emph{Combinatorial Optimization: Polyhedra and Efficiency}, Springer, Berlin, 2003.

\end{thebibliography}
\end{document}